\documentclass[11pt,reqno]{amsart}
\usepackage{verbatim}
\usepackage{amsmath}
\usepackage{amsthm}
\usepackage{amsfonts}
\usepackage{mathrsfs}
\usepackage{amssymb}
\usepackage{amsbsy}
\usepackage{graphicx}
\usepackage{setspace}

\newtheorem{theorem}{Theorem}[section]
\newtheorem{corollary}[theorem]{Corollary}

\newtheorem{remark}[theorem]{Remark}
\newtheorem{proposition}[theorem]{Proposition}
\newcommand{\cj}[1]{\overline{#1}}

\usepackage[left=.8in,right=.8in,top=.7in,bottom=.7in,footnotesep=18pt]{geometry}

\newenvironment{demo}[1]{\renewcommand{\proofname}{#1}\begin{proof}}{\end{proof}}

\title{An expansion for polynomials orthogonal over an analytic Jordan curve}

\author{Erwin Mi\~{n}a-D\'{\i}az}
\address{Indiana-Purdue University Fort Wayne, Department of Mathematical Sciences,
2101 E. Coliseum Blvd.,
Fort Wayne, IN 46805-1499, USA.}

\begin{document}

\maketitle

\begin{abstract}We consider polynomials that are orthogonal over
an analytic Jordan curve $L$ with respect to a positive analytic weight, and show that each such polynomial of sufficiently large degree can be expanded in a series of certain
integral transforms that converges uniformly in the
whole complex plane. This expansion yields, in particular and simultaneously, Szeg\H{o}'s classical strong asymptotic formula and a new integral representation for the polynomials inside $L$. We further exploit such a representation to derive finer asymptotic results for weights having finitely many singularities (all of algebraic type) on a thin neighborhood of the orthogonality curve. Our results are a generalization of those previously obtained in \cite{andrei} for the case of $L$ being the unit circle.
\end{abstract}

\section{Introduction and statements of the results}\label{sec3}

The study of polynomials orthogonal over a closed rectifiable curve of the complex plane was initiated by Szeg\H{o} in \cite{Szego1}, and later continued by Szeg\H{o} himself and such authors as Smirnov, Keldysh, Lavrentiev, Korovkin, Suetin and Geronimous (see \cite{Suetin} for references and an overview of the developments until 1964).
Polynomials orthogonal over several arcs and curves have also been studied, for instance (and without being exhaustive), by Akhiezer \cite{akh1}, \cite{akh2}, Widom \cite{Widom}, Aptekarev \cite{aptekarev}, Peherstorfer and coauthors \cite{Pehers5}, \cite{Pehers4}, \cite{Pehers3}, \cite{Pehers2}, \cite{Pehers1}, and for an orthogonality measure with finitely many point masses outside the curve/arc, by Kaliaguine \cite{Kal2}, \cite{Kal1}.

Among the central questions that are often investigated figure the asymptotic behavior of the orthogonal polynomials and the distribution and location of their zeros. In this regard, the case of a closed curve has the peculiarity (not observed in that of an open arc) that the interior of its polynomial convex hull is non-empty\footnote{A well-known result by Widom \cite{Widom2} asserts that the zeros must accumulate, in the limit, on the polynomial convex hull of the support of the orthogonality measure.}, giving more freedom of  distribution to the zeros of the polynomials, and consequently, making the behavior of the polynomials themselves much less clear. The results of the present work clarify this question to a substantial extent for a single closed curve under analyticity conditions.

Let $L_1$ be an analytic Jordan curve in the complex plane $\mathbb{C}$ and let $h(z)$ be an
analytic function in a neighborhood of $L_1$ such that $h(z)>0$ for all
$z\in L_1$. Using the Gram-Schmidt orthogonalization process, we can form a unique sequence $\left\{p_n(z)\right\}_{n=0}^\infty$ of orthonormal polynomials over $L_1$ with respect to $h(z)$, i.e., satisfying
\begin{equation}\label{eq6}
p_n(z)=\gamma_n z^n+\mathrm{lower\ degree\
terms},\quad \gamma_n>0,
\quad n\geq 0,
\end{equation}
\begin{equation}\label{eq7}
\frac{1}{2\pi}\oint_{L_1}p_n(z)\cj{p_m(z)}h(z)|dz|=\left\{\begin{array}{ll}
                                                                 0,\ &\ n\not=m, \\
                                                                 1,\ &\ n=m.
                                                               \end{array}
\right.
\end{equation}

In what follows, we are concerned with the asymptotic behavior of these polynomials as their degree $n$ becomes large. With this generality, essentially the only known result is Szeg\H{o}'s strong asymptotic formula
\begin{equation}\label{eq20}
p_{n}(z)=\Delta_e(z;h)\sqrt{\phi'(z)}[\phi(z)]^n\left[1+o(1)\right].
\end{equation}

Here $\phi$ is the conformal map of the exterior $\Omega_1$ of $L_1$ onto
the exterior of the unit circle satisfying that
$\phi(\infty)=\infty$, $\phi'(\infty)>0$, $\Delta_e(z;h)$ is the so-called exterior Szeg\H{o} function for the weight $h$, and (\ref{eq20}) holds locally uniformly as $n\to\infty $ on any open set  $\Omega_\rho\supset\overline{\Omega}_1$ that is conformally mapped by $\phi$ onto the exterior of a circle about the origin of radius $\rho <1$, and is
such that $\Delta_e(z;h)$ is analytic on $\Omega_\rho$ (see next subsection for details).

For $h(z)\equiv 1$, this formula was established by Szeg\H{o} in his paper \cite{Szego1} of 1921, while for an arbitrary positive analytic weight, it first appears in Chap. XVI of the first edition of his book \cite{Szego} of 1939. So far as this writer can learn, progress in understanding the asymptotic behavior of $p_n(z)$ at the remaining points of the complex plane, that is, for $z\in \mathbb{C}\setminus \Omega_\rho$, has only been made in the specific case of $L_1$ being the unit circle, the strongest result having being obtained recently in \cite{andrei}. Here the authors use the Riemann-Hilbert approach for the asymptotic analysis of orthogonal polynomials  to derive, for each $p_n(z)$ of sufficiently large degree, a series expansion in terms of certain recursively generated Cauchy transforms. This important result yields at once Szeg\H{o}'s asymptotic formula and an integral representation for $p_n(z)$ inside the unit circle, from which it is possible to distill the precise behavior of the polynomials under additional assumptions on the first singularities encountered by the exterior Szeg\H{o} function. This has been  done in \cite{andrei} for finitely many polar singularities, as well as for two examples of an isolated essential singularity. Earlier related works (e.g., \cite{Sza}) are briefly described in the introduction of \cite{andrei}.

In the present paper we extend the expansion of \cite{andrei} to an arbitrary analytic curve. Our proof is not based on the Riemann-Hilbert method, it is rather direct and in some sense natural, which we believe will lead to applications to other systems of orthogonal polynomials. From the dominant term of the expansion we derive precise asymptotic formulas for $p_n(z)$ in a case where the exterior Szeg\H{o} function has finitely many algebraic singularities in a thin neighborhood of $L_1$. We state our results in Subsections \ref{sub2} and \ref{sub3} below, followed by their proofs in Section \ref{sec2}.

\subsection{Preliminaries} In this subsection we introduce some notation to be used throughout, as well as the concepts involved in the asymptotic behavior of $p_n$. In particular, we discuss the Szeg\H{o} functions associated with the weight $h$. For a deeper discussion of these functions we refer the reader to Chap. X of \cite{Szego}.

Given $r\geq 0$, we set
\[
\mathbb{T}_r:=\{w:|w|=r\},\quad
\mathbb{E}_r:=\{w:r<|w|\leq \infty\},\quad
\mathbb{D}_r:=\{w:|w|<r \}.
\]

If $K$ is a set and $f(z)$ a function defined on
$K$, then $\overline{K}$ and $\partial K$ denote, respectively, the closure and the boundary of $K$, and
$f(K):=\{f(z):z\in K\}$.

\paragraph{Szeg\H{o} functions}

Let $f(t)$ be an analytic function defined on a neighborhood of the unit circle $\mathbb{T}_1$ such that $f(t)>0$ for all $t\in \mathbb{T}_1$. The function
\begin{equation}\label{eq5}
w\mapsto \exp\left\{\frac{1}{4\pi}\oint_{\mathbb{T}_1}\log
f(t)\frac{t+w}{t-w}|dt|\right\}=\exp\left\{\frac{1}{4\pi i}\oint_{\mathbb{T}_1}\frac{\log
f(t)}{t}\cdot\frac{t+w}{t-w}dt\right\},\quad w\in
\cj{\mathbb{C}}\setminus \mathbb{T}_1,
\end{equation}
is analytic on $\cj{\mathbb{C}}\setminus
\mathbb{T}_1$.  Its restriction to
$\mathbb{D}_1$ is called the \emph{interior}
Szeg\H{o} function for $f$, and we denote it by
$D_i(w;f)$. It is univocally determined by the
properties
\begin{enumerate}
\item[(a)] $D_i(w;f)$ has an analytic continuation from $\mathbb{D}_1$ to some neighborhood of $\cj{\mathbb{D}}_1$,
$D_i(w;f)\not=0$ for all $w
\in\cj{\mathbb{D}}_1$ and $D_i(0;f)>0$;

\item[(b)]  $|D_i(w;f)|^{2}=f(w)$ for all $w\in \mathbb{T}_1$.
\end{enumerate}

Property (a) easily follows by noticing that $\log f(t)$ is analytic in a neighborhood of $\mathbb{T}_1$, and therefore, the analytic continuation of $D_i(w;f)$ to $\cj{\mathbb{D}}_1$ is given by the expression to the right of the $=$ sign in (\ref{eq5}) if integration is taken over a circle about the origin of radius slightly larger that $1$. Property (b) is a consequence of that $|D_i(w;f)|^2$ is the exponential of the Poisson integral of $\log f(t)$ (see, e.g., \cite[Thm. 11.8]{Rudin}).

The restriction of the function in (\ref{eq5}) to $\mathbb{E}_1$
is called the \emph{exterior} Szeg\H{o} function
for $f$ and we denote it by $D_e(w;f)$. Notice that
\[
D_e(w;f)=\frac{1}{\cj{D_i(1/\cj{w}\,;f)}},\quad w\in \mathbb{E}_1,
\]
so that $D_e(w;f)$ is univocally determined by the properties
\begin{enumerate}
\item[(a$'$)] $D_e(w;f)$ has an analytic continuation from $\mathbb{E}_1$ to some neighborhood of $\cj{\mathbb{E}}_1$; $D_e(w;f)\not=0$ for all $z\in\cj{\mathbb{E}}_1$ and $D_e(\infty;f)>0$;

\item[(b$'$)] $|D_e(w;f)|^{-2}=f(w)$ for all $w\in \mathbb{T}_1$.
\end{enumerate}

These considerations can be generalized to an arbitrary analytic Jordan curve $L_1$ as follows. Let $\Omega_1$
be the exterior of $L_1$, that is, the unbounded component of $\overline{\mathbb{C}}\setminus
L_1$, and let
$\psi=\psi(w)$ be the unique conformal map of
$\mathbb{E}_1$ onto $\Omega_1$ satisfying that
$\psi(\infty)=\infty$, $\psi'(\infty)>0$.

Let $\widehat{\rho}\geq 0$ be the smallest number such that
$\psi$ has an analytic and \emph{univalent}
continuation from $\mathbb{E}_1$ to
$\mathbb{E}_{\widehat{\rho}}$.  Because $L_1$ is analytic, $\widehat{\rho}<1$.

For every $r$ with $\widehat{\rho}\leq r<\infty$, set
\[
\Omega_r:=\psi(\mathbb{E}_r),\quad L_r:=\partial
\Omega_r,\quad G_r:=\mathbb{C}\setminus
\overline{\Omega}_r,
\]
and let
\[
\phi(z):\Omega_{\widehat{\rho}}\mapsto
\mathbb{E}_{\widehat{\rho}}
\]
be the inverse function of $\psi$. Observe that for every $r>\widehat{\rho}$, $L_r$ is an
analytic Jordan curve.

Then, given a weight function $h(z)$ that is positive and analytic on $L_1$, the \emph{exterior} Szeg\H{o} function $\Delta_e(z;h)$ for $h$ is defined as
\begin{equation*}
\Delta_e(z;h):=D_e(\phi(z);h\circ \psi)=\exp\left\{\frac{1}{4\pi i}\oint_{L_1}\log
h(\zeta)\frac{1+\cj{\phi(\zeta)}\phi(z)}{\phi(\zeta)-\phi(z)}\phi'(\zeta)d\zeta\right\},\quad z\in\Omega_1,
\end{equation*}
which is uniquely determined by the properties
\begin{enumerate}
\item[(i$'$)] $\Delta_e(z;h)$ is analytic and never zero on
$\cj{\Omega}_1$, $\Delta_e(\infty;h)>0$;

\item[(ii$'$)]  $|\Delta_e(z;h)|^{-2}=h(z)$ for all $z\in L_1$.
\end{enumerate}

Similarly, any conformal map $\varphi(z)$ of the interior domain $G_1$ of $L_1$ onto
the unit disk $\mathbb{D}_1$ has an analytic and univalent continuation to some neighborhood of $\cj{G}_1$.  Denoting its inverse by $\delta(w):\cj{\mathbb{D}}_1\mapsto \cj{G}_1$, an \emph{interior} Szeg\H{o} function $\Delta_i(z;h)$ for $h$ is defined as
\begin{equation}\label{eq3}
\Delta_i(z;h):=D_i\left(\varphi(z);h\circ
\delta\right)=\exp\left\{\frac{1}{4\pi i}\oint_{L_1}\log
h(\zeta)\frac{1+\cj{\varphi(\zeta)}\varphi(z)}{\varphi(\zeta)-\varphi(z)}\varphi'(\zeta)d\zeta\right\},\quad  z\in G_1,
\end{equation}
which satisfies the properties
\begin{enumerate}
\item[(i)] $\Delta_i(z;h)$ is analytic and never zero on
$\cj{G}_1$;

\item[(ii)] $|\Delta_i(z;h)|^{2}=h(z)$ for all $z\in L_1$.
\end{enumerate}

Although the interior Szeg\H{o} functions as defined by (\ref{eq3}) depend on the choice of $\varphi$, any two of them differ at most in a multiplicative constant of modulus $1$. Hereafter we shall assume that one such $\Delta_i(z;h)$ has been fixed.

Let $\rho$ be the smallest number larger than or equal to $\widehat{\rho}$ such that $\Delta_e(z;h)$ is analytic on $\Omega_{\rho}$. By property (i$'$), $\widehat{\rho}\leq \rho<1$. Moreover,
\begin{enumerate}
\item[(iii)] $\Delta_i(z;h)^{-1}$ has an analytic continuation from $\cj{G}_1$ to $G_{1/\rho}$.
\end{enumerate}

To see this, consider the Schwarz function $S(z)$ of the curve $L_1$ (see, e.g., \cite{Davis}), which is analytic and univalent in a neighborhood of $L_1$, and is univocally determined by the property that $S(z)=\cj{z}$ for $z\in L_1$. Indeed, $S(z)$ is well-defined all over the band $\Omega_{\widehat{\rho}}\cap G_{1/\widehat{\rho}}$, where it can be expressed in terms of the exterior conformal maps as
\[
S(z)=\cj{\psi\left(1/\cj{\phi(z)}\right)},\quad z\in\Omega_{\widehat{\rho}}\cap G_{1/\widehat{\rho}}.
\]
Let
\[
z^*:=\cj{S(z)}, \quad z\in \Omega_{\widehat{\rho}}\cap G_{1/\widehat{\rho}},
\]
be the so-called Schwarz reflection of $z$ about $L_1$. Then, the analytic continuation of $\Delta_i(z;h)^{-1}$ from $\cj{G}_1$ to $G_{1/\rho}$ is given by
\[
\frac{1}{\Delta_i(z;h)}:=\Delta_e(z;h)\cj{\Delta_e\left(z^*;h\right)}
\cj{\Delta_i\left(z^*;h\right)},\quad z\in
G_{1/\rho}\setminus\cj{G}_1.
\]

\paragraph{The kernel $W(\zeta,z)$} Let $\varphi$ denote, as above, a conformal map of $G_1$ onto $\mathbb{D}_1$, and define the meromorphic kernel
\begin{equation}\label{eq4}
W(\zeta,z):=\frac{\sqrt{\varphi'(z)}\sqrt{\varphi'(\zeta)}}{\varphi(\zeta)-\varphi(z)},\quad
\zeta,\,z\in G_1.
\end{equation}

That this kernel does not depend on the choice of $\varphi$ can be easily verified from the fact that any other conformal
map $\varphi_1$  of $G_1$ onto
$\mathbb{D}_1$ is related to $\varphi$ through a  M\"obius
transformation, that is,
\[
\varphi(z)=e^{i\theta}\frac{\varphi_1(z)-\varphi_1(z_0)}{1-\overline{\varphi_1(z_0)}\varphi_1(z)},\qquad
e^{i\theta}=\varphi'(z_0)\left(1-|\varphi_1(z_0)|^2\right)/\varphi'_1(z_0),
\]
where $z_0$ is that point of $G_1$ mapped by
$\varphi$ onto $0$.

Moreover, if we choose, as we may, a conformal map $\varphi$ that does not vanish on $\Omega_{\widehat{\rho}}\cap G_1$, then this $\varphi$ has an analytic and univalent continuation from $\cj{G}_1$ to $G_{1/\widehat{\rho}}$ given by
\[
\varphi(z)=\frac{1}{\overline{\varphi\left(z^*\right)}},
\qquad z\in G_{1/\widehat{\rho}}\setminus \overline{G}_1,
\]
so that $W(\zeta,z)$ can be extended as a function
\[
W(\zeta,z):G_{1/\widehat{\rho}}\times G_{1/\widehat{\rho}}\to
\overline{\mathbb{C}}
\]
in such a way that for every fixed $z\in G_{1/\widehat{\rho}}$,
$W(\zeta,z)$ is analytic in the variable $\zeta$ on $G_{1/\widehat{\rho}}\setminus
\{z\}$, with a simple pole at $z$ of residue $1$.

We finish this subsection by noticing that positive analytic weights $h(z)$ over $L_1$ are easy to generate because they are precisely those of the form
\[
h(z)=V(z)\cj{V(z^*)},\quad z\in L_1,
\]
with $V(z)$ a zero free analytic function in a neighborhood of $L_1$.

\subsection{The expansions}\label{sub2}

Hereafter we will suppress $h$ and simply write
$\Delta_e(z)$ for $\Delta_e(z;h)$ and
$\Delta_i(z)$ for $\Delta_i(z;h)$.

Fix a number $r$ such that $\rho<r<1$, and for each
integer $n\geq 0$, let us recursively define  the
following sequence of functions:
\[
f_n^{(0)}(z):= 1, \quad z\in \cj{\mathbb{C}},
\]
and for all $k\geq 0$,
\begin{equation*}
f_n^{(2k+1)}(z):= -\frac{1}{2\pi i
}\oint_{L_{r}}f_n^{(2k)}(\zeta)\Delta_e(\zeta)\Delta_i(\zeta)W(\zeta,z)\sqrt{\phi'(\zeta)}
\,[\phi(\zeta)]^nd\zeta,\quad z\in
G_{1/\rho}\setminus L_r,
\end{equation*}
\begin{equation*}
f_n^{(2k+2)}(z):= \frac{1}{2\pi
i}\oint_{L_{1/r}}\frac{f_n^{(2k+1)}(\zeta)\sqrt{\phi'(\zeta)}\,[\phi(\zeta)]^{-n}d\zeta}{\Delta_e(\zeta)\Delta_i(\zeta)[\phi(\zeta)-\phi(z)]},\quad
z\in \Omega_{\rho}\setminus L_{1/r}.
\end{equation*}

Let
\[
\Lambda_r:=\max_{\zeta\in
L_r}\left|\frac{\Delta_e(\zeta)\Delta_i(\zeta)}{\sqrt{\phi'(\zeta)}}\right|,\qquad
\Lambda'_{r}:=\max_{\zeta\in
L_{1/r}}\left|\sqrt{\phi'(\zeta)}\Delta_e(\zeta)\Delta_i(\zeta)\right|^{-1},
\]
and
\[
M_r:=\max_{(\zeta,z)\in L_r\times
L_{1/r}}|W(\zeta,z)|<\infty,
\]
so that obviously (verify it by mathematical induction), for all $k\geq 0$ ,

\begin{equation}\label{eq70}
\quad\left|f^{(2k+1)}_n(z)\right|\leq \Lambda_r
r^{n+1}
\left[\frac{\Lambda_r\Lambda'_rM_rr^{2n}}{1/r-r}\right]^{k}\max_{\zeta\in
L_r}|W(\zeta,z)|,\quad z\in G_{1/\rho}\setminus
L_r,\quad k\geq 0,
\end{equation}
\begin{equation}\label{eq71}
 \left|f^{(2k+2)}_n(z)\right|\leq
\frac{\Lambda_{r}\Lambda'_{r}M_r r^{2n}}{|1/r-|\phi(z)||}
\left[\frac{\Lambda_r\Lambda'_rM_rr^{2n}}{1/r-r}\right]^k,\quad
z\in \Omega_{\rho}\setminus L_{1/r}, \quad
k\geq 0.
\end{equation}

It follows that the two series
\[
f^{(0)}_n(z)+f^{(2)}_n(z)+f^{(4)}_n(z)+\cdots+f^{(2k)}_n(z)+\cdots,\quad
z\in \Omega_{\rho}\setminus L_{1/r},
\]
\[
f^{(1)}_n(z)+f^{(3)}_n(z)+\cdots+f^{(2k+1)}_n(z)+\cdots,\quad
z\in G_{1/\rho}\setminus L_r,
\]
converge absolutely and locally uniformly in
their respective regions of definition, provided
$n$ is so large that
\begin{equation}\label{eq1}
\frac{\Lambda_r\Lambda'_rM_rr^{2n}}{1/r-r}<1.
\end{equation}

Let $P_n(z)$ be the $n$th monic orthogonal polynomial, that is,
\[
P_n(z)=\gamma_n^{-1}p_n(z),\quad n\geq 0,
\]
where $p_n$ satisfies (\ref{eq6}) and (\ref{eq7}).

\begin{theorem}\label{thm1}
For every $n$ satisfying (\ref{eq1}), we have
\begin{equation}\label{eq76}
\frac{\Delta_e(\infty)P_n(z)}{[\phi'(\infty)]^{-n-1/2}}=\left\{\begin{array}{ll}
                   {\displaystyle
\Delta_e(z)\sqrt{\phi'(z)}\,[\phi(z)]^n\sum_{k=0}^\infty
f_n^{(2k)}(z),} &\quad
z\in \Omega_{1/r},\\
{\displaystyle
\Delta_e(z)\sqrt{\phi'(z)}\,[\phi(z)]^n\sum_{k=0}^\infty
f_n^{(2k)}(z)-\frac{1}{\Delta_i(z)}\sum_{k=0}^\infty
f_n^{(2k+1)}(z)}, &\quad z\in \Omega_{r}\cap G_{1/r}\,, \\
{\displaystyle
-\frac{1}{\Delta_i(z)}\sum_{k=0}^\infty
f_n^{(2k+1)}(z)}, & \quad z\in G_{r},
\end{array}
\right.
\end{equation}
and
\begin{equation}\label{eq74}
\gamma_n^{-2}=\Delta_e(\infty)^{-2}[\phi'(\infty)]^{-2n-1}\left(1+\sum_{k=0}^\infty
\frac{1}{2\pi i
}\oint_{L_{1/r}}\frac{f_n^{(2k+1)}(\zeta)\sqrt{\phi'(\zeta)}\,[\phi(\zeta)]^{-n-1}d\zeta}{\Delta_e(\zeta)
\Delta_i(\zeta) }\right).
\end{equation}
\end{theorem}

Let us now consider the following slightly
different sequence of integral transforms. For each fixed integer $n\geq 0$, set
\[
g_n^{(0)}(z):= 1, \quad z\in \cj{\mathbb{C}},
\]
and for all $k\geq 0$,
\begin{eqnarray*}
g_n^{(2k+1)}(z)&:=& -\frac{1}{2\pi i
}\oint_{L_{r}}g_n^{(2k)}(\zeta)\Delta_e(\zeta)\Delta_i(\zeta)
W(\zeta,z)\sqrt{\phi'(\zeta)}\,[\phi(\zeta)]^nd\zeta,\quad
z\in G_{1/\rho}\setminus L_r,\\ &&
\\ g_n^{(2k+2)}(z)&:=&\frac{\phi(z)}{2\pi
i}\oint_{L_{1/r}}\frac{g_n^{(2k+1)}(\zeta)\sqrt{\phi'(\zeta)}\,[\phi(\zeta)]^{-n-1}d\zeta}
{\Delta_e(\zeta)\Delta_i(\zeta)[\phi(\zeta)-\phi(z)]},\quad
z\in \Omega_{\rho}\setminus L_{1/r}.
\end{eqnarray*}
Then,
\begin{equation*}
\quad\left|g^{(2k+1)}_n(z)\right|\leq \Lambda_r
r^{n+1}
\left[\frac{\Lambda_r\Lambda'_rM_rr^{2n+2}}{1/r-r}\right]^{k}\max_{\zeta\in
L_r}|W(\zeta,z)|,\quad z\in G_{1/\rho}\setminus
L_r,\quad k\geq 0,
\end{equation*}
\begin{equation*}
 \left|g^{(2k+2)}_n(z)\right|\leq
\frac{|\phi(z)|\Lambda_{r}\Lambda'_{r}M_r r^{2n+1}}{|1/r-|\phi(z)||}
\left[\frac{\Lambda_r\Lambda'_rM_rr^{2n+2}}{1/r-r}\right]^k,\quad
z\in \Omega_{\rho}\setminus L_{1/r}, \quad
k\geq 0,
\end{equation*}
and the following theorem holds true:

\begin{theorem}\label{thm2}
For every $n$ satisfying (\ref{eq1}), we have
\begin{equation*}
\frac{\gamma_n^2P_n(z)}{\Delta_e(\infty)[\phi'(\infty)]^{n+1/2}}=\left\{\begin{array}{ll}
                   {\displaystyle
\Delta_e(z)\sqrt{\phi'(z)}\,[\phi(z)]^n\sum_{k=0}^\infty
g_n^{(2k)}(z)}, &\
z\in \Omega_{1/r},\\
 {\displaystyle
\Delta_e(z)\sqrt{\phi'(z)}\,[\phi(z)]^n\sum_{k=0}^\infty
g_n^{(2k)}(z)-\frac{1}{\Delta_i(z)}\sum_{k=0}^\infty
g_n^{(2k+1)}(z)}, &\ z\in \Omega_{r}\cap G_{1/r}, \\
 {\displaystyle -\frac{1}{\Delta_i(z)}\sum_{k=0}^\infty
g_n^{(2k+1)}(z)}, & \ z\in G_{r}.
\end{array}
\right.
\end{equation*}
In particular,
\begin{equation}\label{eq69}
\gamma_n^2=
\Delta_e(\infty)^2[\phi'(\infty)]^{2n+1}\sum_{k=0}^\infty
g_n^{(2k)}(\infty).
\end{equation}
\end{theorem}

\begin{remark} These expansions have been previously obtained in \cite{andrei} for
$L_1=\mathbb{T}_1$. They are the outcome of applying the steepest descent method of Deift and Zhou for the asymptotic analysis of a matrix Riemann-Hilbert problem solved by the orthogonal polynomials and closely related functions. Theorem \ref{thm1} extends
Theorem 1 of \cite{andrei}, while relation
(\ref{eq69}) extends Theorem 2 of \cite{andrei}. Our proof of Theorem
\ref{thm1} (and similarly, that of Theorem
\ref{thm2}) is direct: call $H_n(z)$
the function in the right-hand side of
(\ref{eq76}) and observe that the three
expressions that define it in the corresponding
components of $\cj{\mathbb{C}}\setminus
\left(L_r\cup L_{1/r}\right)$ are redundant, in
the sense that they are analytic continuations of
each other. Thus, $H_n(z)$ is an entire function
with a pole of order $n$ at $\infty$, therefore
it is a polynomial. Proving that it is orthogonal
to all powers of $z^m$, $0\leq m<n$, is also
straightforward.
\end{remark}

\begin{corollary}\label{cor1} Let $r_1$ be such that
$\rho<r_1<1/\rho$. Then for every $r>\rho$,
\begin{equation}\label{eq73}
\gamma_n=\Delta_e(\infty)[\phi'(\infty)]^{n+1/2}\left\{1+\mathcal{O}\left(r^{2n}\right)\right\}\,,
\end{equation}
\begin{equation}\label{eq67}
P_n(z)=
\frac{\Delta_e(z)\sqrt{\phi'(z)}\,[\phi(z)]^n}{\Delta_e(\infty)[\phi'(\infty)]^{n+1/2}}\left\{1+\mathcal{O}\left(\frac{r^n}{r_1^n}\right)\right\},
\quad z\in \cj{\Omega}_{r_1}.
\end{equation}

If $\rho<r_1<1$, then
\begin{equation}\label{eq68}
\frac{\Delta_e(\infty)P_n(z)}{[\phi'(\infty)]^{-n-1/2}}=\frac{\Delta_i(z)^{-1}}{2\pi
i
}\oint_{L_{1}}\Delta_e(\zeta)\Delta_i(\zeta)W(\zeta,z)\sqrt{\phi'(\zeta)}\,[\phi(\zeta)]^nd\zeta
+\mathcal{O}\left(r_1^nr^{2n}\right),\quad z\in
\cj{G}_{r_1}\,.
\end{equation}
Equalities (\ref{eq67}) and (\ref{eq68}) hold
uniformly as $n\to\infty$.
\end{corollary}

Of course, (\ref{eq73}) and (\ref{eq67}) are equivalent to
\begin{equation}\label{eq10}
p_n(z)=
\Delta_e(z)\sqrt{\phi'(z)}\,[\phi(z)]^n\left\{1+\mathcal{O}\left(\frac{r^n}{r_1^n}\right)\right\},
\quad z\in \cj{\Omega}_{r_1}.
\end{equation}

Formula (\ref{eq10}) is due to Szeg\H{o} \cite[Thm. 16.5]{Szego}, though he established it with a less precise estimate for the rate of decay of the error term. For $\rho< r_1 \leq 1$, the estimate in (\ref{eq10}) was already obtained by Suetin in \cite{Suetin} (see formula (2.16) in there).

\subsection{Positive weights with algebraic singularities near $L_1$}\label{sub3}

We can derive from (\ref{eq68}) finer asymptotic formulas for $P_n$ if we know more about the  singularities of both the exterior Szeg\H{o} function and the map $\psi$. For instance, if $h(z)\equiv 1$, then $\Delta_e(z)\equiv 1$, $\Delta_i(z)\equiv 1$, and the behavior of $P_n$ inside $L_1$ only depends on geometric considerations and can be determined with great precision, for instance, when $\partial\Omega_{\widehat{\rho}}$ is a piecewise analytic curve, in which case the map $\psi$ has finitely many singularities on the circle $\mathbb{T}_{\widehat{\rho}}$, having an asymptotic expansion about each of them. We will not pursue the analysis of this case here as it is very similar to the one already carried out in \cite{mina} for polynomials orthogonal over the interior of an analytic curve with respect to area measure. Instead, we shall concentrate on a case where the behavior of $P_n$ is only influenced by the singularities of $\Delta_e(z)$, which are finitely many, all lying on the band $G_1\cap\Omega_{\widehat{\rho}}$ and of algebraic type.

Let $a_1,a_2,\ldots,a_s$ be $s\geq 1$
distinct complex numbers all lying on a curve
$L_{\rho}$ with $\widehat{\rho}<\rho<1$. Let
$\lambda_1\geq\lambda_2\geq \cdots\geq\lambda_s$
be such that $\lambda_k\in \mathbb{R}\setminus
\left\{0,-1,-2,\ldots\right\}$ for all $1\leq
k\leq s$, and let $u$ be the number of subindexes $k$ for which $\lambda_k=\lambda_1$, so that
\[
\lambda_1=\lambda_2=\cdots=\lambda_{u}> \lambda_{u+1}\geq \cdots\geq\lambda_s \quad (1\leq u\leq s).
\]

Consider a weight of the form
\begin{equation}\label{eq79}
h(z):=|\omega(z)|^{-2}\prod_{k=1}^s|z-a_k|^{2\lambda_k},\quad
z\in L_1,
\end{equation}
where $\omega(z)$ is an analytic function on
$\cj{\Omega}_{\sigma}$ for some $\widehat{\rho}<\sigma<\rho$, positive at $\infty$ and never zero on
$\cj{\Omega}_{1}\cup \{a_1,a_2,\ldots,a_s\}$.

Let the numbers $\rho_k$, $\sigma_k$, $\Theta_k$ ($1\leq k\leq s$) be defined from the $a_k$'s by the relations
\begin{equation}\label{eq30}
\rho_k:=\phi(a_k)=\rho e^{i\Theta_k},\quad \sigma_k:= \sigma e^{i\Theta_k}, \quad 0\leq
\Theta_k<2\pi,\quad 1\leq k\leq s.
\end{equation}

Let $[\sigma_k,\rho_k]$ be the segment joining $\sigma_k$ and $\rho_k$, and define (see Figure \ref{Fig1} below)
\begin{equation}\label{eq31}
\Gamma_\sigma:=\mathbb{T}_\sigma\cup\{\rho_k:\lambda_k\in\mathbb{N}\}\cup \left(\cup_{ \lambda_k\not \in\mathbb{N}}[\sigma_k,\rho_k]\right),\quad
\Sigma_\sigma:=\left\{z\in\Omega_\sigma: \phi(z)\not\in \Gamma_\sigma\right\},
\end{equation}
so that the exterior Szeg\H{o} function
for the weight $h$ in (\ref{eq79}) is analytic
on $\Sigma_\sigma$ with $a_1,a_2,\ldots,a_s$
being its only singularities on $L_{\rho}$, since, indeed,
\begin{equation}\label{e}
\Delta_e(z)=\omega(z)\prod_{k=1}^s\left(\frac{\phi(z)}{z-a_k}\right)^{\lambda_k},\quad
z\in \Sigma_\sigma,
\end{equation}
with the branches of the $\lambda_k$-power
functions chosen so that
$[\phi'(\infty)]^{\lambda_k}>0$.

\begin{figure}
\centering 
\includegraphics[scale=.6]{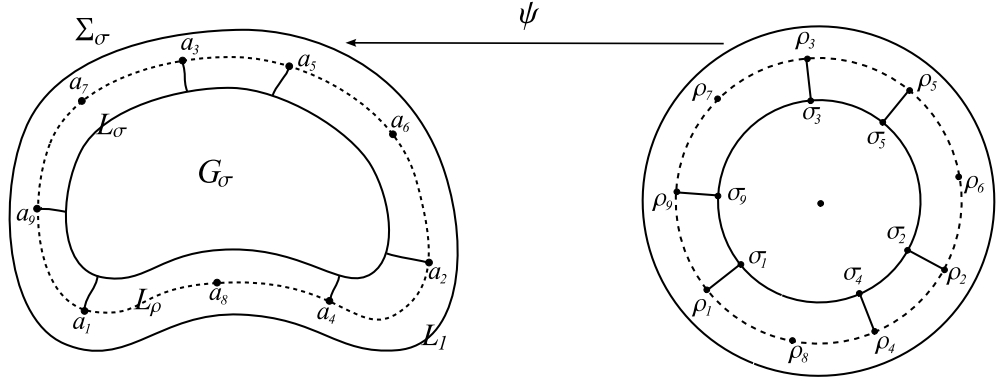}
\caption{Illustration of 9 singularities, $a_6$, $a_7$ and $a_8$ are poles.}\label{Fig1}
\end{figure}

An interior Szeg\H{o} function
$\Delta_i(z)$ for $h$ is given by
\begin{equation}\label{SzegoEq9}
\Delta_i(z)=\Delta_i\left(z;|\omega|^{-2}\right)\prod_{k=1}^s\left(\frac{(z-a_k)\left[1-\cj{\varphi(a_k)}\varphi(z)\right]}
{\varphi(z)-\varphi(a_k)}\right)^{\lambda_k}
\,,\quad
 z\in \cj{G}_{1}\,.
\end{equation}

In what follows,
\[
\alpha_k:=
[\phi'(a_k)]^{\lambda_k-1/2}\lim_{z\to a_k}\left\{\left(\frac{\phi(z)}{z-a_k}\right)^{-\lambda_k}\Delta_e(z)\right\},\quad 1\leq k\leq s,
\]
and $\binom{a}{b}$ stands for the generalized binomial coefficient, i.e.,
\[
\binom{a}{b}:=\frac{\Gamma(a+1)}{\Gamma(b+1)\Gamma(a-b+1)},
\]
where $\Gamma$ denotes the Euler gamma function.

\begin{theorem}\label{thm3} (a) For all $z\in G_1 \setminus \partial\Sigma_\sigma$,
\begin{eqnarray}\label{eq35}
\frac{\Delta_e(\infty)P_n(z)}{[\phi'(\infty)]^{-n-1/2}}&=&
\left\{\begin{array}{ll}
\Delta_e(z)\sqrt{\phi'(z)}[\phi(z)]^n,\ &\  z\in G_1\cap \Sigma_\sigma, \\
0,\ &\ z\in G_\sigma,
\end{array}
\right.\nonumber\\
&&+\binom{n}{\lambda_1-1}\Delta_i(z)^{-1}\left(\sum_{k=1}^u\alpha_k\Delta_i(a_k)
W(a_k,z)[\phi(a_k)]^{n+1}+\rho^nr_n(z)\right)
\end{eqnarray}
where $r_n(z)$ converges locally uniformly to $0$ on $
G_1 \setminus \partial\Sigma_\sigma$.

(b) For every $1\leq j\leq s$,
\begin{eqnarray}\label{eq36}
\frac{\Delta_e(\infty)P_n(a_j)}{[\phi'(\infty)]^{-n-1/2}}&=&
\binom{n}{\lambda_j}\alpha_j\phi'(a_j)\phi(a_j)^{n} +\binom{n}{\lambda_1-1}\Delta_i(a_j)^{-1}\sum_{1\leq k\leq u,\, k\not=j}\alpha_k\Delta_i(a_k)
W(a_k,a_j)[\phi(a_k)]^{n+1}\nonumber\\
&&+o\left(n^{\max\{\lambda_j,\lambda_1-1\}}\rho^n\right)
\end{eqnarray}
as $n\to\infty$.
\end{theorem}

Several remarks are in order.

\begin{remark} The proof of Theorem \ref{thm3} yields the following estimates for the rate of decay of the functions $r_n(z)$ in (\ref{eq35}). If $\lambda_1=1$ and $u=s$, then for every compact set $K\subset G_1 \setminus \partial\Sigma_\sigma$, there is $0<\delta<1$ such that $r_n(z)=\mathcal{O}(\delta^n)$ uniformly on $K$ as $n\to\infty$. Otherwise,
\[
r_n(z)=\left\{\begin{array}{lll}
{\displaystyle \mathcal{O}\left(n^{\lambda_{u+1}-\lambda_1}\right)},\ &\ if\ \lambda_1=1,\ u<s,\  \\
{\displaystyle
\mathcal{O}\left(n^{-1}\right)},\ &\
if\ \lambda_1\not=1,\ u=s,\\
{\displaystyle
\mathcal{O}\left(n^{-\min\left\{1,\,\lambda_1-\lambda_{u+1}\right\}}\right)},\ &\
if\ \lambda_1\not=1,\ u<s,
\end{array}
\right.
\]
locally uniformly on $G_1 \setminus \partial\Sigma_\sigma$ as $n\to\infty$.

Likewise, a better and generally exact estimate for the $o$-error term in (\ref{eq36}) can be easily obtained from  the proof of (\ref{eq36}), though a somewhat tedious case comparison is required.
\end{remark}

\begin{remark} Many fine results on the location and distribution of the zeros of the polynomials $P_n(z)$ follow from Theorem \ref{thm3}(a). For instance:
\begin{enumerate}
\item[(a)] For every compact $K\subset \Omega_\rho$ there is a number $N_K$ such that if $n\geq N_K$, then $P_n(z)$ has no zeros on $K$, and this is also true for any compact $K\subset \cj{\Omega}_\rho\setminus\{a_1,\ldots,a_s\}$ whenever $\lambda_1<1$;
\item[(b)] For every compact $K\subset G_\rho$ there is a number $N_K$ such that if $n\geq N_K$, then $P_n(z)$ has at most $u-1$ zeros on $K$, counting multiplicities.
\item[(c)] Let $\mathcal{Z}$ be the set of those points $t\in\cj{\mathbb{C}}$ such that every neighborhood of $t$ contains zeros of infinitely many polynomials $P_n(z)$. Then, $\cj{\Omega}_\rho\cap \mathcal{Z}=L_\rho$.

\item[(d)] $\mathcal{Z}\cap G_{\rho}$ consists of those points $t\in G_\rho$ satisfying an equation of the form
\[
\sum_{k=1}^u\alpha_k\Delta_i(a_k)
W(a_k,t)e^{i\theta_k}=0,
\]
with angles $\theta_1,\ldots,\theta_u$ for which it is possible to find a subsequence $\{n_j\}_{j\geq 1}\subset\mathbb{N}$ such that
\begin{equation}\label{eq17}
e^{i\theta_k}=\lim_{j\to \infty}e^{i(n_j+1)\Theta_k},\quad 1\leq k\leq u.
\end{equation}

\item[(e)] For each $n\geq 1$, let $\mu_n$ be the normalized counting measure of the zeros $z_{n,1},z_{n,2},\ldots, z_{n,n}$ of $P_n$, that is, $\mu_n:=n^{-1}\sum_{k=1}^n\delta_{z_{n,k}}$, where $\delta_z$ denotes the Dirac unit point measure at $z$. Let $\mu_{L_\rho}$ be the \emph{equilibrium measure} of the compact set $L_\rho$, whose value at any given Borel set
$B\subset L_\rho$ is
\begin{equation}\label{eq37}
\mu_{L_\rho}(B):=\frac{1}{2\pi\rho}\int_{\phi(B)}|dt|.
\end{equation}
Then, the sequence $\{\mu_n\}_{n\geq 1}$ converges in the weak*-topology to $\mu_{L_\rho}$, i.e., for every continuous function
$f$ defined on $\overline{\mathbb{C}}$,
$\lim_{n\to\infty}\int fd\mu_n=\int fd\mu_{L_\rho}$.
\end{enumerate}

Moreover, a result similar to Theorem 4 of \cite{andrei} (see also Theorems 11.1 and 11.2 of \cite{Simon}) on the separation and speed of convergence to $L_\rho$ of the zeros of $P_n$ can be also obtained from (\ref{eq35}).

\end{remark}

Statement (a) actually follows straightforwardly  from (\ref{eq35}) and the maximum modulus principle, while (b) and (d) are also easy if one gets help from Hurwitz theorem. Which $u$-tuples $\{\theta_1, \ldots,\theta_u\}$ satisfy (\ref{eq17}) depends on the specific values of the angles $\Theta_k$ and can be characterized as it has been done in \cite[Thm. 5]{Sza}, \cite[Prop. 3]{andrei} (or, if more details are needed, see \cite[Sec. 2.2]{mina}). Statement (c) is a clear consequence of (e), given that, by definition (\ref{eq37}), the support of $\mu_{L_\rho}$ is $L_\rho$. The proof of (e) is based on standard arguments of logarithmic potential theory: by (a), (b) and (\ref{eq67}), any measure $\mu$ that is the weak*-limit of some subsequence $\{\mu_{n_j}\}_{j\geq 1}$ is supported on $L_\rho$ and satisfies
\[
\int\log\frac{1}{|z-t|}d\mu(t)=\int\lim_{j\to\infty}\log\frac{1}{|z-t|}d\mu_{n_j}(t)=\lim_{j\to\infty}{n_j}^{-1}\log |P_{n_j}(z)|^{-1}=\log|\phi'(\infty)/\phi(z)|,\quad z\in \Omega_\rho.
\]
On the other hand, it is not difficult to verify from (\ref{eq37}) that
\[
\int\log\frac{1}{|z-t|}d\mu_{L_\rho}(t)=\log|\phi'(\infty)/\phi(z)|,\quad z\in \Omega_\rho,
\]
i.e., the  logarithmic potential of $\mu$ coincides outside $L_\rho$ with that of $\mu_{L_\rho}$, which, by a well-known theorem of Carleson \cite[Thm. 4.13]{ST}, implies that $\mu=\mu_{L_\rho}$.

\begin{remark}Values of
$\lambda_k\in\left\{0,-1,-2,\ldots\right\}$ are
purposely excluded because their corresponding
factors $(z-a_k)^{\lambda_k}$ would not create a
singularity (but a zero) for $\Delta_e(z)$ at $a_k$,
and therefore, these factors may be simply
regarded as being part of the function $h(z)$. We also
note that among the weights defined by (\ref{eq79}) are those of the form
\begin{equation*}
h(z):=|\omega_1(z)|^{-2}\prod_{k=1}^s|z-a_k|^{2\lambda_k},\quad
 \{a_1,a_2,\ldots,a_s\}\subset
L_{\rho}\cup L_{1/\rho},
\end{equation*}
with $\omega_1(z)$ an analytic function on
$\cj{\Omega}_{\rho}$ that is never zero on
$\cj{\Omega}_{1}\cup \{a_k\in
L_{\rho}\}\cup\{a_k^*:a_k\in L_{1/\rho} \}$ (here $a^*_k$ denotes the Schwarz reflection of $a_k$ about $L_1$). For in such a case, $h(z)$ can be also written in
the form (\ref{eq79}) as follows:
\begin{equation*}
h(z)=|\omega(z)|^{-2}\left(\prod_{k\,:\,a_k\in
L_{\rho}}|z-a_k|^{2\lambda_k}\right)\left(\prod_{k\,:\,a_k\in
L_{1/\rho}}|z-a^*_k|^{2\lambda_k}\right),\quad
z\in L_1,
\end{equation*}
with
\[
\omega(z)=\omega_1(z)\prod_{k\,:\,a_k\in
L_{1/\rho}}|\phi(a_k)|^{\lambda_k}\left(\frac{\phi(z)-\phi(a_k^*)}{z-a^*_k}
\cdot\frac{z-a_k}{\phi(z)-\phi(a_k)}\right)^{\lambda_k}\,.
\]
\end{remark}

\section{Proofs}\label{sec2}

\begin{demo}{Proof of Theorem \ref{thm1}} Let us denote by $H_n(z)$ the right-hand side of
(\ref{eq76}), which is originally defined on
$\cj{\mathbb{C}}\setminus\left(L_{r}\cup
L_{1/r}\right)$, and let us prove that $H_n(z)$
is indeed an entire function.

Let
\begin{equation}\label{eq12}
H^+_n(z):=\Delta_e(z)\sqrt{\phi'(z)}\,[\phi(z)]^n\left[1+\sum_{k=0}^\infty
\frac{1}{2\pi
i}\oint_{L_{1}}\frac{f_n^{(2k+1)}(\zeta)\sqrt{\phi'(\zeta)}\,[\phi(\zeta)]^{-n}d\zeta}
{\Delta_e(\zeta)\Delta_i(\zeta)[\phi(\zeta)-\phi(z)]}\right], \quad z\in \Omega_1,
\end{equation}
which, in view of (\ref{eq70}), is
well-defined and analytic on $\Omega_1$. This function provides the analytic continuation of
$\left.H_n\right|_{\Omega_{1/r}}$ to $\Omega_1$, which  follows from the very definition of $\left.H_n\right|_{\Omega_{1/r}}$, given that for all $k\geq 0$ and $z\in \Omega_{1/r}$ (and by deforming
$L_{1/r}$ into $L_1$),
\begin{eqnarray*}
f^{(2k+2)}_n(z):=\frac{1}{2\pi
i}\oint_{L_{1/r}}\frac{f_n^{(2k+1)}(\zeta)\sqrt{\phi'(\zeta)}\,[\phi(\zeta)]^{-n}d\zeta}
{\Delta_e(\zeta)\Delta_i(\zeta)[\phi(\zeta)-\phi(z)]}=\frac{1}{2\pi
i}\oint_{L_{1}}\frac{f_n^{(2k+1)}(\zeta)\sqrt{\phi'(\zeta)}\,[\phi(\zeta)]^{-n}d\zeta}
{\Delta_e(\zeta)\Delta_i(\zeta)[\phi(\zeta)-\phi(z)]}.
\end{eqnarray*}

Moreover, by the residue theorem (deforming $L_1$ back into $L_{1/r}$ in (\ref{eq12})), we find that for
every $z\in\Omega_{1}\cap G_{1/r}$,
\begin{equation*}
H_n^+(z)=
\Delta_e(z)\sqrt{\phi'(z)}\,[\phi(z)]^n\left[1+\sum_{k=0}^\infty
\frac{1}{2\pi
i}\oint_{L_{1/r}}\frac{f_n^{(2k+1)}(\zeta)\sqrt{\phi'(\zeta)}\,[\phi(\zeta)]^{-n}d\zeta}
{\Delta_e(\zeta)\Delta_i(\zeta)[\phi(\zeta)-\phi(z)]}\right]-\frac{1}{\Delta_i(z)}\sum_{k=0}^\infty
f_n^{(2k+1)}(z),
\end{equation*}
that is, the analytic continuation $H^+_n$ of
$\left.H_n\right|_{\Omega_{1/r}}$ to $\Omega_1$
coincides for values of
$z\in\Omega_{1}\cap G_{1/r}$ with $\left.H_n\right|_{\Omega_{r}\cap G_{1/r}}$  as defined by
(\ref{eq76}).

Similarly,
\begin{equation*}
H_n^-(z):=\frac{1}{\Delta_i(z)}
\sum_{k=0}^\infty\frac{1}{2\pi
i }\oint_{L_1}f_n^{(2k)}(\zeta)\Delta_e(\zeta)\Delta_i(\zeta)W(\zeta,z)
\sqrt{\phi'(\zeta)}\,[\phi(\zeta)]^nd\zeta,\quad z\in G_1,
\end{equation*}
provides the analytic continuation of $\left.H_n\right|_{G_r}$ to $G_1$, which for values of $z\in \Omega_r \cap G_1$ coincides precisely with $\left.H_n\right|_{\Omega_{r}\cap G_{1/r}}$.

Thus, $H_n(z)$ is an
entire function and
\begin{eqnarray*}
\lim_{z\to\infty}\frac{H_n(z)}{z^n}&=&
\Delta_e(\infty)[\phi'(\infty)]^{n+1/2}\left(1+f^{(2)}_n(\infty)+f^{(4)}_n(\infty)+\cdots+f^{(2k)}_n(\infty)+\cdots \right)\\
&=&\Delta_e(\infty)[\phi'(\infty)]^{n+1/2}.
\end{eqnarray*}
By Liouville's theorem, $H_n(z)$ is a polynomial
of exact degree $n$, whose leading coefficient is
\[
\Delta_e(\infty)[\phi'(\infty)]^{n+1/2}.
\]

Now, from the definition of $\left.H_n\right|_{\Omega_{r}\cap G_{1/r}}$ we have
\begin{eqnarray*}
\frac{1}{2\pi}\oint_{L_1}H_n(z)\cj{z^m}h(z)|dz|&=&
\sum_{k=0}^\infty\frac{1}{2\pi}\oint_{L_1}\sqrt{\phi'(z)}\,[\phi(z)]^n\cj{z^{m}}f^{(2k)}_n(z)\Delta_e(z)h(z)|dz|\\
&&-\sum_{k=0}^\infty\frac{1}{2\pi}\oint_{L_1}f^{(2k+1)}_n(z)\cj{z^{m}}\Delta_i(z)^{-1}h(z)|dz|,
\end{eqnarray*}
so that Theorem \ref{thm1} will follow at once if we show that
\begin{eqnarray}\label{eq66}
&&\frac{1}{2\pi}\oint_{L_1}\sqrt{\phi'(z)}\,[\phi(z)]^n\cj{z^{m}}f^{(2k)}_n(z)\Delta_e(z)h(z)|dz|
\nonumber\\
&=&\left\{\begin{array}{ll}
0,\ &\ 0\leq m<n,\ k\geq 0, \\
{\displaystyle
\Delta_e(\infty)^{-1}[\phi'(\infty)]^{-n-1/2}},\
&\  m=n,\ k=0,\\
            {\displaystyle
\frac{\Delta_e(\infty)^{-1}[\phi'(\infty)]^{-n-1/2}}{2\pi i
}\oint_{L_{1/r}}\frac{f_n^{(2k-1)}(\zeta)\sqrt{\phi'(\zeta)}[\phi(\zeta)]^{-n-1}d\zeta}{\Delta_e(\zeta)\Delta_i(\zeta)
}},\ &\
            m=n,\ k\geq 1,
          \end{array}
\right.
\end{eqnarray}
and
\begin{equation}\label{eq15}
\frac{1}{2\pi}\oint_{L_1}f^{(2k+1)}_n(z)\cj{z^{m}}\Delta_i(z)^{-1}h(z)|dz|
=0,\quad n,m\geq 0,\quad k\geq 0.
\end{equation}

First, we obtain by making the change of variables
$z=\psi(w)$ that for all $0\leq m\leq n$,
\begin{equation}\label{eq65}
\frac{1}{2\pi}\oint_{L_1}\sqrt{\phi'(z)}\,[\phi(z)]^n\cj{z^{m}}f^{(2k)}_n(z)\Delta_e(z)h(z)|dz|\\
=\frac{1}{2\pi}\oint_{\mathbb{T}_1}\cj{\left[\frac{\sqrt{\psi'(w)}[\psi(w)]^{m}}{\Delta_e(\psi(w))}\right]}w^nf^{(2k)}_n(\psi(w))
|dw|,
\end{equation}
where we have used that for $z\in L_1$, $h(z)=|\Delta_e(z)|^{-2}$.

Now, the function $\sqrt{\psi'(w)}[\psi(w)]^{m}/\Delta_e(\psi(w))$
is analytic on
$\cj{\mathbb{E}}_{1}\setminus\{\infty\}$ with a
pole of order $m$ at $\infty$, so that from its Laurent expansion at infinity we obtain that for certain coefficients $a_j$ (that  depend on $n$ and $m$),
\begin{equation}\label{eq63}
\cj{\left[\frac{\sqrt{\psi'(w)}[\psi(w)]^{m}}{\Delta_e(\psi(w))}\right]}w^n
=\Delta_e(\infty)^{-1}[\phi'(\infty)]^{-m-1/2}w^{n-m}\left(1+\sum_{j=1}^\infty
a_jw^{j}\right),\quad w\in\mathbb{T}_1.
\end{equation}

On the other hand, from the definition of $f^{(2k)}_n(z)$ for $k\geq 1$ we see that
\begin{eqnarray*}
f^{(2k)}_n(\psi(w))&=&\frac{1}{2\pi
i}\oint_{L_{1/r}}\frac{f_n^{(2k-1)}(\zeta)\sqrt{\phi'(\zeta)}[\phi(\zeta)]^{-n}d\zeta}
{\Delta_e(\zeta)\Delta_i(\zeta)
[\phi(\zeta)-w]}\\
&=&\frac{1}{2\pi
i}\oint_{\mathbb{T}_{1/r}}\frac{f_n^{(2k-1)}(\psi(t))t^{-n}\sqrt{\psi'(t)}dt}
{\Delta_e(\psi(t))\Delta_i(\psi(t))
(t-w)}, \quad w\in
\mathbb{E}_{\rho}\setminus \mathbb{T}_{1/r},
\end{eqnarray*}
is indeed analytic in all of
$\cj{\mathbb{C}}\setminus\mathbb{T}_{1/r}$, and we obtain from its Taylor expansion about $0$ that for certain coefficients $b_j$ (that  depend on $n$ and $k$),
\begin{equation}\label{eq64}
 f^{(2k)}_n(\psi(w))=\frac{1}{2\pi
i
}\oint_{L_{1/r}}\frac{f_n^{(2k-1)}(\zeta)\sqrt{\phi'(\zeta)}[\phi(\zeta)]^{-n-1}d\zeta}{\Delta_e(\zeta)
\Delta_i(\zeta)}+\sum_{j=1}^\infty b_jw^j\,,\quad
 w\in \mathbb{T}_1.
\end{equation}
Taking into account that $f^{(0)}_n(z)\equiv 1$, we then get (\ref{eq66}) by combining (\ref{eq65}),
(\ref{eq63}) and (\ref{eq64}).

Similarly, if $\varphi$ is a conformal map of
$G_1$ onto $\mathbb{D}_1$ and $\delta(w)$ is its
inverse, we have
\begin{equation}\label{eq16}
\frac{1}{2\pi}\oint_{L_1}f^{(2k+1)}_n(z)\cj{z^{m}}\Delta_i(z)^{-1}h(z)|dz|
=\frac{1}{2\pi}\oint_{\mathbb{T}_1}f^{(2k+1)}_n(\delta(w))\sqrt{\delta'(w)}\,\cj{\sqrt{\delta'(w)}
[\delta(w)]^{\,m}
\Delta_i(\delta(w))
}|dw|,
\end{equation}
where we have used that for $z\in L_1$, $h(z)=|\Delta_i(z)|^{2}$.

On the one hand, $\sqrt{\delta'(w)}[\delta(w)]^{m}\Delta_i(\delta(w))$ is
analytic on $\cj{\mathbb{D}}_1$, and from its Taylor expansion we obtain
\begin{equation}\label{eq14}
\cj{\sqrt{\delta'(w)}[\delta(w)]^{m}\Delta_i(\delta(w))}=\sum_{j=0}^{\infty}c_j
w^{-j},\quad w\in \mathbb{T}_1.
\end{equation}

On the other hand, from the definition of $f^{(2k+1)}_n(z)$ and (\ref{eq4}) we have that
\[
f^{(2k+1)}_n(\delta(w))\sqrt{\delta'(w)}=-\frac{1}{2\pi
i}\oint_{L_{r}}\frac{f_n^{(2k)}(\zeta)\Delta_e(\zeta)\Delta_i(\zeta)\sqrt{\varphi'(\zeta)}
\sqrt{\phi'(\zeta)}\,[\phi(\zeta)]^nd\zeta}{\varphi(\zeta)-w}
\quad
\]
is analytic on $\cj{\mathbb{C}}\setminus\varphi(L_r)\supset
\cj{\mathbb{E}}_1$ and vanishes at $\infty$, so that its Laurent expansion at $\infty$ restricted to $\mathbb{T}_1$ is of the form
\begin{equation}\label{eq13}
f^{(2k+1)}_n(\delta(w))\sqrt{\delta'(w)}=\sum_{j=1}^\infty d_j w^{-j},\quad
w\in \mathbb{T}_1.
\end{equation}

Thus, (\ref{eq15}) follows by inserting (\ref{eq14}) and (\ref{eq13}) in  (\ref{eq16}).
\end{demo}

\begin{demo}{Proof of Theorem \ref{thm2}} Proceed just as in the proof of Theorem \ref{thm1} above.
\end{demo}

\begin{demo}{Proof of Corollary \ref{cor1}} Let $r_1$ and
$r$ be such that $\rho<r_1<1/\rho$,  $\rho<r<\min\left\{r_1,r_1^{-1}\right\}$, so that
$L_{r_1}\subset G_{1/r}\cap\Omega_{r}$.

By inequalities (\ref{eq70}) and (\ref{eq71}), we
see that
\begin{equation*}
\sum_{k=1}^\infty\left|f^{(2k+1)}_n(z)\right|\leq
r^{3n+1}\frac{\Lambda_r^2\Lambda_r'M_r\max_{(\zeta,z)\in
L_r\times
L_{r_1}}|W(\zeta,z)|}{1/r-r-\Lambda_r\Lambda_r'M_rr^{2n}},\quad
z\in L_{r_1},
\end{equation*}
\begin{equation*}
\sum_{k=0}^\infty
\left|f^{(2k+2)}_n(z)\right|\leq r^{2n}
\frac{(1/r-r)\Lambda_{r}\Lambda'_{r}M_r }{(1/r-r_1)(1/r-r-\Lambda_r\Lambda_r'M_rr^{2n})}
,\quad
 z\in L_{r_1},
\end{equation*}
and we obtain from Theorem \ref{thm1} that
\begin{equation}\label{eq72}
\frac{\Delta_e(\infty)P_n(z)}{[\phi'(\infty)]^{-n-1/2}}=
\Delta_e(z)\sqrt{\phi'(z)}\,[\phi(z)]^n-\Delta_i(z)^{-1}
f_n^{(1)}(z)+\mathcal{O}\left(r_1^nr^{2n}\right),
\quad z\in L_{r_1}.
\end{equation}
Given that, again by (\ref{eq70}),
$f_n^{(1)}(z)=\mathcal{O}\left(r^{n}\right)$ uniformly in
$z\in L_{r_1}$ as $n\to\infty$, we get from (\ref{eq72}) that
(\ref{eq67}) holds uniformly in $z\in L_{r_1}$, and by the
maximum modulus principle for analytic functions,
it also holds on $\cj{\Omega}_{r_1}$.

If now $r_1<1$, then from the definition of $f_n^{(1)}(z)$ and the
residue theorem we obtain that for all $z\in
L_{r_1}$,
\begin{eqnarray*}
f_n^{(1)}(z)&=&-\frac{\Delta_i(z)^{-1}}{2\pi
i
}\oint_{L_{r}}\Delta_e(\zeta)\Delta_i(\zeta)W(\zeta,z)\sqrt{\phi'(\zeta)}
\,[\phi(\zeta)]^nd\zeta\\
&=&
\Delta_e(z)\sqrt{\phi'(z)}\,[\phi(z)]^n-\frac{\Delta_i(z)^{-1}}{2\pi
i
}\oint_{L_{1}}\Delta_e(\zeta)\Delta_i(\zeta)W(\zeta,z)\sqrt{\phi'(\zeta)}
\,[\phi(\zeta)]^nd\zeta,
\end{eqnarray*}
which together with (\ref{eq72}) yields that
(\ref{eq68}) holds for $z\in L_{r_1}$, and again by the
maximum modulus principle for analytic functions,
it also holds on $\cj{G}_{r_1}$.

Equality (\ref{eq73}) follows, for instance, from (\ref{eq74}) and (\ref{eq70}), or from (\ref{eq69}).
\end{demo}

\begin{demo}{Proof of Theorem \ref{thm3}}  We first prove a proposition that will help the proof of Theorem \ref{thm3} to go through smoothly. The following notation will be used. For each $\delta>0$ and $t\in \mathbb{C}$,
\[
D_\delta(t):=\{w:|w-t|<\delta\},\quad T_\delta(t):=\{w:|w-t|=\delta\}=\partial D_\delta(t).
\]

Let $0<\sigma<\rho$ be given numbers, and define $\delta:=\rho-\sigma$. Suppose that  $v(t,z)$ is a function of two complex variables that is analytic in the variable $t$ on the closed disk $\cj{D_{2\delta}(\rho)}$ for each $z\in E$ ($E$ certain set), and that
\[
\sup\left\{|v(t,z)|:(t,z)\in\cj{D_{2\delta}(\rho)}\times E\right\}<\infty,
\]
so that, by the Cauchy integral formula, we also have that for every integer $p\geq 0$, there is a constant $0<M_p<\infty$ such that
\[
|\partial^p
v(t,z)/\partial t^p|\leq M_p,\quad (t,z)\in\cj{D_{\delta}(\rho)}\times E.
\]

For $\beta\in \mathbb{R}\setminus
\left\{0,-1,-2,\ldots\right\}$, let the function $(t-\rho)^{-\beta}$ be defined for $t\in \mathbb{C}\setminus (-\infty,\rho]$ according to the branch of the argument
\[
-\pi<\arg(t-\rho)<\pi\,,\quad
t\in \mathbb{C}\setminus (-\infty,\rho],
\]
and let
\[
 (t-\rho)_-^{-\beta}:=\lim_{z\to t,\, \Im{z}>0}(z-\rho)^{-\beta}, \quad (t-\rho)_+^{-\beta}:=\lim_{z\to t,\, \Im{z}<0}(z-\rho)^{-\beta},\quad  t\in \mathbb{C}\setminus (-\infty,\rho],
\]
denote its boundary values from either side of the semi-line  $(-\infty,\rho]$.

\begin{proposition} With the notations above, and for every integer $n\geq 0$, we have
\begin{equation}\label{eq28}
I :=\frac{1}{2\pi i
}\oint_{T_\delta(\rho)}(t-\rho)^{-\beta} t^nv(t,z)dt=\binom{n}{\beta-1}v(\rho,z)\rho^{n-\beta+1}+\left\{\begin{array}{ll}
                      0, &\ if\ \beta=1, \\
                      \mathcal{O}\left(n^{\beta-2}{\rho}^n\right), &\ if\
                      \beta\not= 1,
                      \end{array}\right.
\end{equation}
uniformly in $z\in E$ as $n\to\infty$.
\end{proposition}

To prove the proposition, we first observe that for every integer
$p\geq 0$,
\begin{eqnarray}\label{eq22}
\frac{\partial^p\left[v(t,z)t^{n}\right]}{\partial
t^p}&=&\sum_{j=0}^p\binom{p}{j}
\frac{\partial^{p-j}v(t,z)}{\partial
t^{p-j}}t^{n-j}\prod_{v=1}^j(n+1-v)=\frac{n!\,t^{n-p}
v(t,z)}{\Gamma(n+1-p)}+\left\{\begin{array}{ll}
                      0, &\ \mathrm{if}\ p=0, \\
                      \mathcal{O}\left(n^{p-1}t^n\right), &\ \mathrm{if}\
                      p>0,
                    \end{array}
\right.\nonumber\\
&=&\mathcal{O}(n^{p}t^n)\quad
\end{eqnarray}
uniformly in $z\in E$ as $n\to\infty$.

Therefore, if $\beta$ is a positive integer,
then we obtain from the Cauchy integral formula and (\ref{eq22}) that
\begin{eqnarray}
I &=& \frac{1}{(\beta-1)!}
\left.\frac{\partial^{\beta-1}\left[v(t,z)t^{n}\right]}{\partial
t^{\beta-1}}\right|_{t=\rho}=\binom{n}{\beta-1}v(\rho,z)\rho^{n-\beta+1}+\left\{\begin{array}{ll}
                      0, &\ \mathrm{if}\ \beta=1, \\
                      \mathcal{O}\left(n^{\beta-2}{\rho}^n\right), &\ \mathrm{if}\
                      \beta>1,
                    \end{array}
\right.
\end{eqnarray}
uniformly in $z\in E$ as $n\to\infty$.

Next, consider a $\beta$ that is not an
integer. Let $\bar{\beta}$ be the smallest \emph{nonnegative} integer not less than $\beta$.
Consecutive integrations by parts over $T_\delta(\rho)$ yield
\begin{eqnarray}\label{eq21}
I&=&\frac{1}{2\pi i}\sum_{j=1}^{\bar{\beta}}\left.\frac{\partial^{j-1}}{\partial
t^{j-1}}\left[v(t,z)t^{n}\right]\right|_{t=\sigma} \frac{(-1)^{j-1}}
{\prod_{l=1}^j(l-\beta)}\left[(\sigma-\rho)_-^{j-\beta}-(\sigma-\rho)_+^{-\beta+j}\right]
\nonumber\\
&&+\frac{(-1)^{\bar{\beta}}}
{2\pi i\prod_{l=1}^{\bar{\beta}}(l-\beta)}\oint_{T_\delta(\rho)}(t-\rho)^{\bar{\beta}-\beta}
\frac{\partial^{\bar{\beta}}}{\partial
t^{\bar{\beta}}}\left[v(t,z)t^{n}\right] dt\,.
\end{eqnarray}

We can now deform $T_\delta(\rho)$ into the two-sided segment $[\sigma,\rho]$ without altering the value of this last integral, and so obtain from (\ref{eq21}) and (\ref{eq22})  that
\begin{eqnarray}\label{eq26}
I&=&
\frac{1}
{2\pi i\prod_{l=1}^{\bar{\beta}}(\beta-l)}\int_{\sigma}^{\rho}
\left[(t-\rho)_+^{\bar{\beta}}-(t-\rho)_-^{\bar{\beta}-\beta}\right]
\frac{\partial^{\bar{\beta}}\left[v(t,z)t^{n}\right]}{\partial
t^{\bar{\beta}}} dt +\sum_{j=1}^{\bar{\beta}}\mathcal{O}\left(n^{j-1}\sigma^n\right)\nonumber\\
&=&\frac{\sin(\pi(\beta-\bar{\beta}))n!\rho^{n-\beta+1}}
{\pi\Gamma(n-\bar{\beta}+1)\prod_{l=1}^{\bar{\beta}}(\beta-l)}\int_{\sigma/\rho}^{1}
(1-x)^{\bar{\beta}-\beta}x^{n-\bar{\beta}}
\left[v(\rho x,z)+
\mathcal{O}\left(n^{-1}\right)\right] dx +\mathcal{O}
\left(n^{\bar{\beta}-1}\sigma^n\right).
\end{eqnarray}

On the one hand, there is some constant $M_1>0$ such that
\begin{equation}\label{eq24}
|v(t,z)-v(\rho,z)|\leq\left|\int_{[\rho,t]}\frac{\partial
v(w,z)}{\partial w}dw\right|\leq M_1|t-\rho|, \quad (t,z)\in\cj{D_{\delta}(\rho)}\times E.
\end{equation}

On the other hand, for every integer $n\geq 0$ and real $\alpha>-1$, we have
\begin{eqnarray*}
\int_{\sigma/\rho}^1(1-x)^\alpha x^{n}
dx&=&\int_0^1(1-x)^\alpha x^{n}
dx-\int_0^{\sigma/\rho}(1-x)^\alpha x^{n}
dx\\
&=&\frac{\Gamma(\alpha+1)\Gamma(n+1)}{\Gamma(n+\alpha+2)}+\mathcal{O}\left(n^{-1}(\sigma/\rho)^{n}\right)=
\frac{\Gamma(\alpha+1)(1+o(1))}{n^{\alpha+1}}\qquad
(n\to\infty),
\end{eqnarray*}
so that by (\ref{eq24}),
\begin{eqnarray}\label{eq25}
&&\int_{\sigma/\rho}^{1}
(1-x)^{\bar{\beta}-\beta}x^{n-\bar{\beta}}
\left[v(\rho x,z)+
\mathcal{O}\left(n^{-1}\right)\right] dx \nonumber\\
&=&v(\rho,z)\int_{\sigma/\rho}^{1}
(1-x)^{\bar{\beta}-\beta}x^{n-\bar{\beta}}dx +n^{-1}\int_{\sigma/\rho}^{1}\mathcal{O}\left((1-x)^{\bar{\beta}-\beta}x^{n-\bar{\beta}}\right)dx\nonumber
+\int_{\sigma/\rho}^{1}\mathcal{O}\left((1-x)^{1+\bar{\beta}-\beta}x^{n-\bar{\beta}}\right)dx\nonumber\\
&=&\frac{\Gamma(\bar{\beta}-\beta+1)\Gamma(n-\bar{\beta}+1)v(\rho,z)}{\Gamma(n-\beta+2)}
+\mathcal{O}\left(n^{-\bar{\beta}+\beta-2}\right) \qquad
(n\to\infty).
\end{eqnarray}

Then, from (\ref{eq26}) and (\ref{eq25}), and taking into account that
\[
\Gamma(\beta-\bar{\beta})\Gamma(\bar{\beta}-\beta+1)=\frac{\pi}{\sin(\pi(\beta-
\bar{\beta}))},\quad  \Gamma(\beta-\bar{\beta})\prod_{l=1}^{\bar{\beta}}(\beta-l)=\Gamma(\beta),
\]
we obtain that (\ref{eq28}) also holds if $\beta$ is not an
integer.

Having proven the proposition above, it is now easy to prove Theorem \ref{thm3}. Let us start by fixing a number $\sigma$ with $\widehat{\rho}<\sigma<\rho$ and such that $\omega(z)$ is analytic on $\cj{\Omega}_\sigma$. Let the corresponding points $\rho_k$, $\sigma_k$ ($1\leq k\leq s$) and the sets $\Gamma_\sigma$ and $\Sigma_\sigma$ be defined as in (\ref{eq30})-(\ref{eq31}).

Let $E\subset G_1\cap \Sigma_\sigma$ and $F\subset G_\sigma$ be compact sets, and let $\rho<r_1<1$ be such that $E\subset \cj{G}_{r_1}$. Then, according to (\ref{eq68}), we have that for all $z\in
\cj{G}_{r_1}\supset E\cup F\cup\{a_1,\ldots,a_s\}$,
\begin{equation}\label{eq80}
\frac{\Delta_e(\infty)P_n(z)}{[\phi'(\infty)]^{-n-1/2}}=\frac{\Delta_i(z)^{-1}}{\,2\pi i
}\oint_{\mathbb{T}_{1}}\Delta_e(\psi(t))\Delta_i(\psi(w))W(\psi(w),z)\sqrt{\psi'(w)}\,w^ndw+
\mathcal{O}\left(r_1^n {\rho}^{3n/2}\right).
\end{equation}

Choose $\widehat{\sigma}$ such that $\sigma<\widehat{\sigma}<\rho$, and if $\delta:=\rho-\widehat{\sigma}$, then the closed disks $\cj{D_{2\delta}(\rho_k)}$, $1\leq k\leq s$, are pairwise disjoint, and
\[
\cj{D_{2\delta}(\rho_k)}\subset \mathbb{D}_1\cap \mathbb{E}_\sigma\setminus\{\phi(z):z\in E\},\quad 1\leq k\leq s.
\]

Define
\[
\widehat{\sigma}_k:=\widehat{\sigma}e^{i\Theta_k}, \quad 1\leq k\leq s,
\]
and the positively oriented contour
\[
\mathcal{C}_{\widehat{\sigma}}:=\mathbb{T}_\sigma\cup
\left(\cup_{\lambda_k\not\in\mathbb{N}}[\sigma_k,\widehat{\sigma}_k]\right)
\cup\left( \cup_{k=1}^s T_\delta(\rho_k)\right),
\]
where each segment $[\sigma_k,\widehat{\sigma}_k]$ is viewed as having two sides.

Since $\omega(z)$ is analytic on $\cj{\Omega}_\sigma$ and
\begin{equation}\label{eq18}
\Delta_e(\psi(w))=\omega(\psi(w))\prod_{k=1}^s\left(\frac{w}{w-\rho_k}\right)^{\lambda_k}\prod_{k=1}^s\left(\frac{w-\rho_k}{\psi(w)-\psi(\rho_k)}\right)^{\lambda_k},\quad
w\in \mathbb{E}_\sigma\setminus \Gamma_\sigma,
\end{equation}
we have that for all $z\in E\cup F\cup \{a_1,a_2,\ldots,a_s\}$, the function (in the variable $w$)
\begin{equation}\label{eq27}
F(w,z):=\Delta_e(\psi(w))\Delta_i(\psi(w))W(\psi(w),z)\sqrt{\psi'(w)}
\end{equation}
is analytic on
\[
\{w:\sigma<|w|<1\}\setminus \Gamma_\sigma
\] (with the exception of the point $\phi(z)$ in case $z\in
E$, where it has a simple pole) with continuous boundary values on
$\mathbb{T}_1\cup \Gamma_\sigma\setminus\{\rho_1,\rho_2,\ldots,\rho_s\}$ when viewing each segment
$[\sigma_k,\rho_k)$ as having two sides. Consequently, by deforming in (\ref{eq80}) $\mathbb{T}_1$ into $\mathcal{C}_{\widehat{\sigma}}$ (and applying the residue theorem in case $z\in E$) we obtain that
\begin{eqnarray}\label{eq19}
&&\frac{\Delta_e(\infty)P_n(z)}{[\phi'(\infty)]^{-n-1/2}}-\left\{\begin{array}{ll}
\Delta_e(z)\sqrt{\phi'(z)}[\phi(z)]^n,\ &\  z\in E, \\
0,\ &\ z\in F\cup \{a_1,a_2,\ldots,a_s\},
\end{array}
\right.\nonumber\\
&=&\frac{\Delta_i(z)^{-1}}{2\pi i
}\oint_{\mathcal{C}_{\widehat{\sigma}}}F(w,z)w^ndw +\mathcal{O}\left(r_1^n {\rho}^{3n/2}\right),\nonumber\\
&=&
\Delta_i(z)^{-1}\sum_{k=1}^s \frac{1}{2\pi i }\oint_{T_\delta(\rho_k)}F(w,z)w^ndw +\mathcal{O}\left({\widehat{\sigma}}^n\right)+\mathcal{O}\left(r_1^n {\rho}^{3n/2}\right),\quad z\in E\cup F\cup \{a_1,\ldots,a_s\}.
\end{eqnarray}

If we now specify
\begin{equation}\label{eq32}
-\pi<\arg(t-\rho)<\pi\,,\quad
t\in \mathbb{C}\setminus (-\infty,\rho]\,,
\end{equation}
\begin{equation}\label{eq33}
-\pi<\arg(t)<\pi\,,\quad t\in
\cj{D_{2\delta}(\rho)},
\end{equation}
we see from (\ref{eq18}), (\ref{eq27}) and  (\ref{eq4}) that for every $1\leq k\leq s$,
\[
F(te^{i\Theta_k},z)=e^{-i\lambda_k\Theta_k}\left(t-\rho\right)^{-\lambda_k}F_k(t,z), \quad t\in \cj{D_{\delta}(\rho)}\setminus [\widehat{\sigma},\rho], \quad z\in E\cup F\cup \left\{a_j: j\not=k\right\}
\]
where $F_k(t,z)$ is analytic (as a function of $t$) on the closed disk $\cj{D_{2\delta}(\rho)}$ for every $z\in E\cup F\cup \left\{\rho_j:j\not=k\right\}$, and
\begin{equation}\label{eq29}
F_k(\rho,z)=\rho_k^{\lambda_k}[\phi'(a_k)]^{\lambda_k-1/2}\Delta_i(a_k)
W(a_k,z)\lim_{z\to a_k}\left\{\left(\frac{\phi(z)}{z-a_k}\right)^{-\lambda_k}\Delta_e(z)\right\}.
\end{equation}
Hence we get from the proposition proven above that for every $z\in E\cup F\cup \left\{a_j:j\not=k\right\}$,
\begin{eqnarray}\label{eq23}
 \frac{1}{2\pi i }\oint_{T_\delta(\rho_k)}F(w,z)w^ndw &=&\frac{e^{i(n-\lambda_k+1)\Theta_k}}{2\pi i }\oint_{T_\delta(\rho)}\left(t-\rho\right)^{-\lambda_k}t^nF_k(t,z)dt\nonumber\\
 &=&\binom{n}{\lambda_k-1}F_k(\rho,z)\rho_k^{n-\lambda_k+1}+\left\{\begin{array}{ll}
                      0, &\ \mathrm{if}\ \lambda_k=1, \\
                      \mathcal{O}\left(n^{\lambda_k-2}{\rho}^n\right), &\ \mathrm{if}\
                      \lambda_k\not= 1.
                      \end{array}\right.
\end{eqnarray}
Thus, part (a) of Theorem \ref{thm3} follows from (\ref{eq19}), (\ref{eq29}) and (\ref{eq23}).

Part (b) follows similarly. The function $W(\psi(w),a_j)$ is analytic in $\{w:\sigma<|w|<1,\ w\not= \rho_j\}$, with a simple pole at $\rho_j$ of residue $\phi'(a_j)$, and with the specifications (\ref{eq32}) and (\ref{eq33}), we have
\[
F(te^{i\Theta_j},a_j)=e^{-i(\lambda_j+1)\Theta_j}\left(t-\rho\right)^{-\lambda_j-1}U_j(t), \quad t\in \cj{D_{\delta}(\rho)}\setminus [\widehat{\sigma},\rho],
\]
where $U_j(t)$ is analytic on $\cj{D_{2\delta}(\rho)}$ and
\begin{equation}\label{eq38}
U_j(\rho)=\rho_j^{\lambda_j}[\phi'(a_j)]^{\lambda_j+1/2}\Delta_i(a_j)
\lim_{z\to a_j}\left\{\left(\frac{\phi(z)}{z-a_j}\right)^{-\lambda_j}\Delta_e(z)\right\},
\end{equation}
so that
\begin{eqnarray}\label{eq34}
 \frac{1}{2\pi i }\oint_{T_\delta(\rho_j)}F(w,a_j)w^ndw &=&\frac{e^{i(n-\lambda_j)\Theta_j}}{2\pi i }\oint_{T_\delta(\rho)}\left(t-\rho\right)^{-\lambda_j-1}t^nU_j(t)dt\nonumber\\
 &=&\binom{n}{\lambda_j}U_j(\rho)\rho_j^{n-\lambda_j}+                      \mathcal{O}\left(n^{\lambda_j-1}{\rho}^n\right),
\end{eqnarray}
and part (b) of the theorem follows by combining (\ref{eq19}), (\ref{eq23}), (\ref{eq38}) and (\ref{eq34}).
\end{demo}


\begin{thebibliography}{99}

\bibitem{akh1}{Akhiezer, N. I.: Orthogonal polynomials on several intervals. Soviet Math. Dokl. \textbf{1},
989-992 (1960)}

\bibitem{akh2}{Akhiezer, N. I. and Tomchuk, Yu. Ya.: On the theory of orthogonal polynomials over several
intervals (in Russian). Dokl. Akad. Nauk SSSR \textbf{138}, 743-746 (1961)}

\bibitem{aptekarev}{Aptekarev, A. I.: Asymptotic properties of polynomials orthogonal
on a system of contours and periodic motions of
Toda chains.  Mat. Sb. \textbf{125}, 231-258 (1984)}

\bibitem{Davis}{Davis, P. J.: \emph{The Schwarz function and its
applications.} The Carus Mathematical Monographs Vol. \textbf{17}. The Mathematical Association of America, 1974.}



\bibitem{Kal2}{Kaliaguine, V. A.: On asymptotics of $L^p$ extremal polynomials on a complex curve $(0<p<\infty)$. J. Approx. Theory \textbf{74}, 226-236 (1993)}

\bibitem{Kal1}{Kaliaguine, V. A.: A note on the asymptotics of orthogonal polynomials on a complex arc: the case of a measure with a discrete part. J. Approx. Theory \textbf{80}, 138-145  (1995)}

\bibitem{andrei}{Mart\'{\i}nez-Finkelshtein, A., McLaughlin, K. T.-R. and Saff, E. B.: Szeg\H{o}
orthogonal polynomials with respect to an
analytic weight: canonical representation and
strong asymptotics. Constr. Approx. \textbf{24},
319-363 (2006)}

\bibitem{mina}{Mi\~{n}a-D\'{\i}az, E.: An asymptotic integral representation for Carleman orthogonal polynomials. Submitted for publication (available at http://arxiv.org/abs/0710.2856)}



\bibitem{Pehers5}{Peherstorfer, F.: On Bernstein-Szeg\H{o} orthogonal polynomials on several intervals. SIAM J.
Math. Anal. \textbf{21}, 461-482  (1990)}

\bibitem{Pehers4}{Peherstorfer, F.: Zeros of polynomials orthogonal on several intervals. Int. Math. Res. Not. \textbf{7}, 361-385 (2003)}

\bibitem{Pehers3}{Peherstorfer, F. and  Yuditskii, P.: Asymptotic behavior of polynomials orthonormal on a homogeneous set. J. Anal. Math. \textbf{89}, 113-154 (2003)}

\bibitem{Pehers2}{Peherstorfer, F.: On the zeros of orthogonal polynomials: the elliptic case. Constr. Approx. \textbf{20}, 377-397 (2004)}

\bibitem{Pehers1}{Lukashov, A. L. and Peherstorfer, F.: Zeros of polynomials orthogonal on two arcs of the unit circle. J. Approx. Theory \textbf{132}, 42-71  (2005)}



\bibitem{Rudin}{Rudin, W.: \emph{Real and complex analysis.} McGraw-Hill, New York, 3rd ed., 1986}

\bibitem{ST}{Saff, E. B. and Totik, V.: \emph{Logarithmic potentials with external
fields.} Berlin: Springer-Verlag, 1997}


\bibitem{Simon}{Simon, B.: Fine structure of the zeros of orthogonal polynomials, I. A tale of two pictures.
 \emph{ETNA} \textbf{25}, 328-368 (2006)}

\bibitem{Suetin}{Suetin, P. K.: Fundamental properties of polynomials orthogonal on a contour. Russ. Math.
Surv. \textbf{21}, 35-83 (1966)}


\bibitem{Sza}{Szabados, J.: On some problems connected with polynomials orthogonal on the complex unit
circle. Act. Math. Scien. Hung. \textbf{33},
197-210 (1979)}

\bibitem{Szego}{Szeg\H{o}, G.: \emph{Orthogonal Polynomials.} Amer. Math. Soc. Colloq. Publ. Vol.
\textbf{23}, Amer. Math. Soc., Providence, RI, 4th ed., 1975}

\bibitem{Szego1}{Szeg\H{o}, G.: \"{U}ber orthogonale polynome, die zu einer gegebenen Kurve der Komplexen Ebene geh\"{o}ren. 	Math. Z. \textbf{9}, 218-270 (1921)}

\bibitem{Widom}{Widom, H: Extremal polynomials associated with a system of curves in the complex plane. Adv. Math. \textbf{3}, 127-232 (1969)}

\bibitem{Widom2}{Widom, H: Polynomials associated with measures in the complex plane. J. Math. Mech. \textbf{16}, 997-1013 (1967)}

\end{thebibliography}
\end{document}